\documentclass[10pt]{article}

\usepackage{amsmath}
\usepackage{amsthm}
\usepackage{amssymb}
\usepackage{amsfonts}
\usepackage{graphicx}
\usepackage{verbatim}
\usepackage{mathrsfs}
\usepackage{subfigure}
\usepackage{fancyhdr}
\usepackage{latexsym}
\usepackage{graphicx}
\usepackage{dsfont}

\theoremstyle{plain}
\newtheorem{theorem}{Theorem}[section]
\newtheorem{proposition}[theorem]{Proposition}
\newtheorem{lemma}[theorem]{Lemma}
\newtheorem{corollary}[theorem]{Corollary}

\newtheorem{remark}[theorem]{Remark}

\theoremstyle{definition}

\newcommand{\al}{\alpha}
\newcommand{\R}{\mathbb{R}}

\newcommand{\E}{\mathcal{E}}
\newcommand{\mR}{\mathcal {R}}

\newcommand{\n}{\noindent}

\newcommand{\bydef}{:=} 

\title{A study of the apsidal angle \\
and a proof of  monotonicity in the logarithmic potential case}
\author{Roberto Castelli \thanks{BCAM - Basque Center for Applied Mathematics, Alameda de Mazarredo, 14
48009 Bilbao, Basque Country, Spain. Email: {\tt rcastelli@bcamath.org}.}
}

\date{}

\begin{document}
\maketitle

\begin{abstract}
This paper concerns the behaviour of the apsidal angle for orbits of central force system with homogenous potential of degree $-2\leq \al\leq 1$ and logarithmic potential. We derive a formula for the apsidal angle as a fixed-end points integral and we study the derivative of the apsidal angle with respect to the angular momentum $\ell$. The monotonicity of the apsidal angle as function of $\ell$ is discussed and it  is proved in the logarithmic potential case.
\end{abstract}
\begin{center}
{\bf \small Keywords} \\ \vspace{.05cm}
{ \small Two-body problem $\cdot$ Logarithmic potential $\cdot$ Apsidal angle  }
\end{center}

\begin{center}
{\bf \small Mathematics Subject Classification (2010)} \\ \vspace{.05cm}
{ \small 70F05 $\cdot$ 65G30 }
\end{center}
\section{Introduction}

In this paper we investigate  the apsidal precession for  the orbits of central force systems of the form
\begin{equation}\label{eq:system}
\ddot u=\nabla V_{\al}(|u|),\quad u\in\R^{2}
\end{equation}
where
\begin{equation}
\begin{array}{ll}
{\displaystyle V_{\al}(x)=\frac{1}{\al}\frac{1}{x^{\al}},\quad }& \al\in[-2,1]\setminus \{0\}\\
\\
{\displaystyle V_{0}(x)=-\log(x) }&\al=0.
\end{array}
\end{equation}
The solutions of system \eqref{eq:system} preserve the mechanical energy $E=\tfrac{1}{2}|\dot u|^{2}-V$ and the angular momentum $\vec\ell=u\wedge \dot u$ and, while rotating,  all the bounded non-collision orbits oscillates between  the apses, that is, the points of minimal ( pericenter)  and maximal (apocenter ) distance from the origin.  The apsidal angle is defined as  the angle at the origin swept by the orbit between two consecutive  apses.

For an admissible choice of $E$ and $\ell$ ( being $\ell$ the scalar angular momentum, positive for counterclockwise  orbits and negative for clockwise orbits), the apsidal angle of an orbit $u(t)$ is given by
\begin{equation}\label{eq:intro}
\Delta_{\al}\theta(u)=\int_{r_{-}}^{r_{+}}\frac{\ell}{r^{2}\sqrt{2\big(E+V(r)\big)-\frac{\ell^{2}}{r^{2}}}}\, dr
\end{equation}
where $r_{\pm}$ are the points where the denominator vanishes.

The extremal cases $\alpha=-2$ and $\alpha=1$ consist in the harmonic oscillator and the Kepler problem respectively.  Systems with $-2<\al<1$ are of relevant interest in quantum mechanics and  astrophysics, for instance in describing  the galactic dynamics in presence of power law densities or massive black-holes and in modelling  the gravitational lensing, see \cite{1997MNRAS, MNR:MNR22071}. Also, the logarithmic potential appears in particle physics \cite{Quigg, 5125922} and in a model for the dynamics of vortex filaments in ideal fluid, see i.e. \cite{Newton}.

Since Newton, the behaviour of the apsidal angle has been  extensively studied, in particular for its implications to celestial mechanics. In the Book I of the {\it Principia} Newton derived a formula that relates the apsidal angle to the magnitude of the attracting force: for ``close to circular'' orbits, a force fields of the form $\kappa r^{n-3}$ leads to an apsidal angle equal to $\pi/\sqrt{n}$. Hence, the experimental measurement of the precession of an orbit may give the exponent of the force law. Newton itself in Book III, looking at the orbit of the Earth and of the Moon, concludes that the attracting force of the Sun and of the Earth must be inverse square of the distance.  ``Close to circular'' means that this formula has been proved for orbits with small eccentricity, or equivalently, with angular momentum close to $\ell_{max}$.

In the singular case ($\al\geq 0$), the behaviour of the apsidal angle plays a fundamental role in dealing with the collision solutions of \eqref{eq:system} and related systems: a possible regularisation of the singularities of the flow is possible only in case that twice the apsidal angle tends to a multiple of $\pi$ as $\ell\to 0$, see for instance \cite{Castelli2013ly, CastelliDCDS, gbegfugfg, Mcgehee}. Also a  variational approach to system \eqref{eq:system} may lead to collision avoidance provided $\Delta_{\al}\theta(u)>\pi$ as $|u|\to 0$, see   \cite{Terracini2007ve}. Moreover, as pointed out in \cite{MNR:MNR8819}, the derivative of the apsidal angle with respect to the angular momentum has common feature with the derivative of the potential scattering phase shift and applications in a variety  of fields, such as nuclear physics and astrophysics.    

 Two are the main technical hurdles that make both  the analytical treatment and the numerical investigation of  the integral \eqref{eq:intro} a difficult problem: the integrand is singular at the end-points and the end-points themselves are only implicitly defined as function of  $E$ and $\ell$.
Nevetherless, partial and important results are very well known.
In 1873, Bertrand proved that there are only two central force laws for which all bounded orbits are closed, namely the linear and the inverse square \cite{Bertrand}. 
It means that the apsidal angle of any bounded solution of the elastic and Kepler problem is rationally proportional to $\pi$ for any value of $E$ and $\ell$. In the first case it values $\pi/2$ and $\pi$ in the second. 
In addition, for $\al=-\tfrac{2}{3}, -1, \tfrac{1}{2}, \tfrac{2}{3}$, the apsidal angle can be expressed in terms of elliptic functions, see \cite{whitta}, but in general $\Delta_{\al}\theta(u)$ can not be written in closed form.

Valluri et al. \cite{MNR:MNR8819} propose an expansion of $\Delta_{\al}\theta(u)$  and  $\partial_{\ell}\Delta_{\al}\theta(u)$ in terms of the eccentricity of the orbit  for   $\alpha$ close $1$.  A study  for the logarithmic case using the $p-ellipse$ approximation and the Lambert W function  is presented in \cite{MNR:MNR22071}. Touma and Tremaine in \cite{1997MNRAS}, by means of the Mellin transform,  provide an asymptotic series on $\ell$ for any $\alpha\in[-2,1]$. Also, by applying a generalisation of the B$\ddot{\rm u}$rmann-Lagrange series for multivalued inverse function, Santos et al. in \cite{PhysRevE.79.036605} obtain a series expression for the apsidal angle  valid for any central force field.

Several numerical simulations, like the ones proposed in the cited papers, show that for any fixed $\al$ and $E$ the apsidal angle is a monotonic function of the angular momentum. However,  at the best of the author's knowledge, no proof of this statement is available in the literature. The purpose of this paper is to review the notion of the apsidal angle and to derive a formula for its derivative with respect to the angular momentum. From where, the monotonicity  will be  evident for any $\al\in(-2,1)$ and proved in the logarithmic potential case.

The plan of the paper is the following: in the next section we recall the orbital structure of the solutions of \eqref{eq:system} and we write a formula for the apsidal angle as a fixed-ends  integral \eqref{eq:deltatheta}, where the parameter $q$ plays the role of the angular momentum. Then, section \ref{sec:Eal} is dedicated to the analytical study of the integrand, in particular of the function $\E_{\al}(s,q)$.  In section \ref{sec:limits} the results so far  obtained are used to compute the limit values of the apsidal angle  for radial and close to circular orbits. In section \ref{sec:derivative} we compute the derivative of the apsidal angle and we provide insight into the monotonicity in the $\alpha$-homogneous case.  Finally, in section \ref{sec:proof} the apsidal angle for $\al=0$  is proved to be monotonically increasing as function of the angular momentum. In Appendix some technical lemmas are collected.

\section{Orbital structure and apsidal angle}\label{sec:aps-ang}

\begin{figure}[htbp]
\begin{center}
\includegraphics[width=\textwidth]{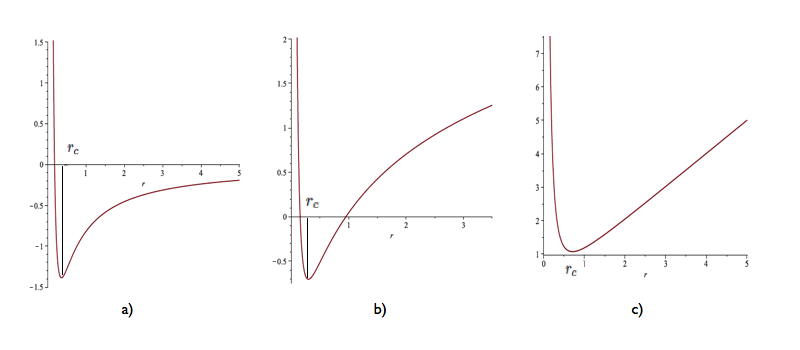}
\caption{{\small Plot of the function $V^{eff}_{\al}(r)$ for $a)\  \al\in(0,1]$, $b)\  \al=0$, $c) \ \al\in[-2,0)$. In each graphic $r_{c}$ is the radius where the minimum of $V_{\al}^{eff}(r)$ is achieved and it corresponds to the radius of the circular orbit. It holds : in $a)$ and $b)$ $r_{c}=\ell^{\frac{2}{2-\al}}$, $V_{\al}^{eff}(r_{c})=-\frac{2-\al}{2\al}\ell^{-\frac{2\al}{2-\al}}$, b) $r_{c}=\ell$, $V_{0}^{eft}(r_{c)}= \frac{1}{2}+\log(\ell).$}}
\label{fig:Veff}
\end{center}
\end{figure}

In polar coordinates $u=(r,\theta)$ the energy associated to a solution $u(t)$ of the problem $\eqref{eq:system}$ is
$$
E(u)=\frac{1}{2}\dot r^{2}+\frac{1}{2}\frac{\ell^{2}}{r^{2}}-V_{\al}(r)
$$
hence the motion is allowed only for those $r$ where $E>V_{\al}^{eff}=\frac{1}{2}\frac{\ell^{2}}{r^{2}}-V_{\al}(r)$.
From the analysis of the $V_{\al}^{eff}$, see also Fig. \ref{fig:Veff}, we realise that
\begin{itemize}
\item{}for $\al\in(0,1]$: bounded orbits exist for negative energies and, for any value $E<0$, the modulus of the  angular momentum $\ell $ ranges between $0$ and $\ell_{max}=\left(\frac{2\al}{\al-2}E\right)^{-\frac{2-\al}{2\al}}$. Values $E\geq 0$ give rise to unbounded orbits;
\item{}for $\al=0$: since $\log(r)\to \infty $ as $r\to \infty$, all the possible orbits are bounded. For any value of $E$, the motion is possible provided $0\leq |\ell|\leq \ell_{max}=e^{E-\frac{1}{2}}$;
\item{}for $\al\in[-2,0)$: the potential function $V_{\al}$ is negative therefore the motion exists only for positive values of $E$. In particular all the orbits are bounded and the $|\ell|$ ranges within the interval $[0,\ell_{\max}]$ with $\ell_{max}=\left(\frac{2\al}{\al-2}E\right)^{-\frac{2-\al}{2\al}}$.
\end{itemize}

The radial coordinate $r(t)$ of the orbit $u(t)$  solves the differential equation
$$
\ddot r-\frac{\ell^{2}}{r^{3}}+\frac{1}{r^{\al+1}}=0
$$ 
and oscillates periodically between the extremal   values $r_{+}$ and $r_{-}$ given as positive solutions of  the apses  equation $ \frac{1}{2}\frac{\ell^{2}}{r^{2}}-V_{\al}(r)-E=0$. While oscillating, the solution rotates around the centre of attraction, in direction according to the sign of $\ell$.

From the relation $\ell=r^{2}\dot \theta$ it follows that the rate of rotation between two consecutive apses is given by $\int_{r_{-}}^{r_{+}}\frac{\ell}{r^{2}\dot r}\, dr$, from where formula \eqref{eq:intro} comes.
\begin{remark}
The apsidal angle $\Delta_{\al}\theta(u)$ is odd as function of $\ell$.
\end{remark}
\begin{remark}
From  \eqref{eq:intro} it descends that $\Delta_{\beta}\theta(u)=\frac{2-\al}{2}\Delta_{\al}\theta(u)$, where $\beta=-\frac{2\al}{2-\al}$, see also \cite{grant}.  Thus apsidal angle for $\al\in[-2,0)$ is determined by the value for $\al\in(0,1]$. 
\end{remark}
Inlight of the two remarks and since the case $\alpha=1$ is completely solved, henceforth we restrict our analysis to
$$
\al\in[0,1),\quad \ell\in(0,1).
$$

A standard approach when studying the apsidal angle is to consider the differential equation 
$$
\frac{d^{2}}{d\theta^{2}}z+z=\phi\left(\frac{1}{z}\right)\frac{1}{z^{2}\ell^{2}}
$$
where $z=\frac{1}{r}$ and $\phi(r)$ is the central force. From here, following  Bertrand \cite{Bertrand} and Griffin \cite{griffin}, it holds

\begin{equation}\label{eq:grif}
\Delta_\alpha\theta(u)=\int_a^b\frac{\sqrt{w(b)-w(a)}}
{\sqrt{a^2w(b)-b^2w(a)+(b^2-a^2)w(z)-z^2(w(b)-w(a))}}dz
\end{equation}
where the function $w(z)=2\int\psi(z)\, dz$ and $\psi(z)=r^{2}\phi(r)$. 
In practice, for the force field we are interest in, it holds
$$
\phi(r)=-\frac{1}{r^{\alpha+1}}, \quad \psi(z)=z^{\alpha-1}, \quad w(z)=\frac{2}{\alpha}z^{\alpha},\qquad \alpha\neq 0
$$
$$
\phi(r)=-\frac{1}{r}, \quad \psi(z)=z^{-1}, \quad w(z)=2\log(z),\qquad \alpha=0.
$$
The dependence on $\ell$ and on the energy of the integral \eqref{eq:grif}  is through the values of the $a=1/r_{+}$ and  $b=1/r_{-}$. Indeed, in the logarithmic potential case, such values are given as solutions of the equation  $\log(z^{2})-\ell^{2}z^{2}=-2E$, whose roots can only be 
  expressed as
$$
z=\sqrt{\frac{-W_{j}(-\ell^{2}e^{-2E})}{\ell}}
$$
where $W_{j}$ represent the two branches of the Lambert W function, \cite{LambertW}. 
 Similarly, in the homogenous case, $a$ and $b$ are the solutions of $-\ell^{2}z^{2}+\frac{2}{\alpha}z^{\alpha}=-2E$.

From \eqref{eq:grif}, by easy manipulation it follows
\begin{equation}
 \begin{split}
 \Delta_\alpha\theta(u)&=
  \int_a^b\frac{dz}{\sqrt{b^{2}-z^{2}-(b^{2}-a^{2})\frac{w(b)-w(z)}{w(b)-w(a)}}}= \int_a^b\frac{dz}{\sqrt{(b-z)\left((b+z)- (b+a)\frac{\frac{w(b)-w(z)}{b-z}}{\frac{w(b)-w(a)}{b-a}}\right)}}\nonumber\\
&=\int_a^b\frac{1}{\sqrt{z-a}}\frac{1}{\sqrt{b-z}}\frac{1}{\sqrt{1+\varepsilon(z)}}\, dz\ \label{integdelta}
\end{split}
\end{equation}
where
\begin{equation}\label{ep}
\varepsilon(z)=\frac{b+a}{z-a}\left[ 1- \frac{\frac{w(b)-w(z)}{b-z}}{\frac{w(b)-w(a)}{b-a}}\right]\ .
\end{equation}
Let $L\bydef b-a$ and perform the change of variable $s=\frac{z-a}{L}$; it results
\begin{equation}\label{intdeltaL}
\Delta_{\alpha}\theta(u)=\int_0^1\frac{1}{\sqrt{s(1-s)}}\frac{1}{\sqrt{1+\epsilon(s)}}\, ds
\end{equation}
where
$$
\epsilon(s)=\frac{(2b-L)}{Ls}\left[1-\frac{w(b)-w\big(b-L(1-s)\big)}{(1-s)\big(w(b)-w(b-L)\big)}\right]\ .
$$
Then, taking $q=\frac{L}{b}$, we obtain

\begin{equation}\label{eq:deltatheta}
\Delta_{\alpha}\theta(u)=\int_0^1\frac{1}{\sqrt{s(1-s)}}\frac{1}{\sqrt{1+\mathcal E_{\alpha}(s,q)}}ds
\end{equation}
with
\begin{equation}\label{eps}
\begin{array}{ll}
	{\displaystyle \mathcal E_{\alpha}(s,q)=\frac{2-q}{qs}\left[ 1-\frac{1-\big(1-q(1-s)\big)^{\alpha}}{(1-s)(1-(1-q)^{\alpha})}\right], }&\alpha\in(0,1)\\
	\\
	{\displaystyle\E_{0}(s,q)=\frac{2-q}{qs}\left[1-\frac{\log(1-q(1-s))}{(1-s)\log(1-q)} \right],\quad }&\alpha=0
\end{array}
\end{equation}

For any fixed value of the mechanical energy, the parameter $q=1-\frac{a}{b}$ is a function of the angular momentum $\ell$. In particular, as $\ell $ increases  from zero to $\ell_{\max}$, the ratio $a/b$ increases from $0$ to $1$. Therefore
$$
q=q(\ell)\in(0,1),\qquad \frac{dq}{d\ell}<0, \quad \forall\ell\in(0,1).
$$  
The integral in \eqref{eq:deltatheta} is the main object of investigation of the present paper.  Two are the reasons why it is advantageous to analyse \eqref{eq:deltatheta} instead of the classical one \eqref{eq:intro}: in the first integral  the end-points  are fixed  and the  dependence on the angular momentum is confined  and well separated from the singularities of the integrand function.

The behaviour of the apsidal angle as the angular momentum $\ell$ varies is intimately related to the behaviour of the function $\E_{\al}(q,s)$ when $s$, $q$ vary in their domain. Hence, we first  present a detailed analysis of the  function $\E_{\al}(s,q)$ and its  derivative $\partial_{q}\E_{\al}(s,q)$.

\section{Analysis of the functions $\E_{\alpha}(s,q)$}\label{sec:Eal}

 It is convenient to introduce the weights
\begin{equation}\label{eq:w}
\omega^{\alpha}_{n}\bydef\left\{\begin{array}{cl}
-{\al \choose n}(-1)^{n} & \al\in(0,1)\\
\\
\frac{1}{n} &\al=0
\end{array}\right.\qquad \forall n\geq 1
\end{equation}
and two auxiliary functions: for any $n\geq 2$ let 
$$
A_{n}(s)\bydef \frac{1-(1-s)^{n-1}}{s},\quad s\in(0,1)
$$
and, for any $n\geq 2$ and $\al\in [0,1)$ let 
\begin{equation}\label{eq:Kn}
\begin{split}
K^{\al}_{n}(s)&\bydef 2\frac{\omega^{\al}_{n+1}}{\omega^{\al}_{n}}A_{n+1}(s)-A_{n}(s)=2\frac{n-\al}{n+1}A_{n+1}(s)-A_{n}(s), \quad s\in(0,1).
\end{split}
\end{equation}
Note that $\omega^{\al}_{n}>0$ for any $n\geq 1$ and $\al\in[0,1)$. We collect in the  Appendix all the properties of the function $A_{n}(s)$ and $K_{n}^{\al}(s)$  used throughout  the paper.

Recalling the power series expansions:
$$
(1-x)^{\al}=1+\sum_{n\geq 1}{\al \choose n}(-1)^{n}x^{n}, \quad \log(1-x)=-\sum_{n\geq 1}\frac{x^{n}}{n}, \quad |x|<1
$$
it holds, for $\al\in(0,1)$,
\begin{equation}
\begin{split}
\mathcal E_{\alpha}(s,q)&=\frac{2-q}{qs}\left[ \frac{-(1-s)\sum_{n\geq 1}{\al \choose n}(-1)^{n}q^{n}+\sum_{n\geq 1}{\al \choose n}(-1)^{n}q^{n}(1-s)^{n}}{(1-s)\left(-\sum_{n\geq 1}{\al \choose n}(-1)^{n}q^{n}\right)}\right]\\
&=\frac{2-q}{qs\left(-\sum_{n\geq 1}{\al \choose n}(-1)^{n}q^{n}\right)}\left[ -\sum_{n\geq 2}{\al \choose n}(-1)^{n}q^{n}+\sum_{n\geq 2}{\al \choose n}(-1)^{n}q^{n}(1-s)^{n-1}\right]\\
&= \frac{1}{\sum_{n\geq 1}{\al \choose n}(-1)^{n}q^{n-1}}\sum_{n\geq 2}{\al \choose n}(-1)^{n}q^{n-2}(2-q)A_{n}(s)
\end{split}
\end{equation}
and similarly, for $\al=0$, 
$$
\E_{0}(s,q)=\frac{1}{\sum_{n\geq 1}\frac{1}{n}q^{n-1}}\sum_{n\geq 2}\frac{1}{n}(2-q)q^{n-2}A_{n}(s).
$$
By means of the weights $\omega_{n}^{\al}$,  for any $\al\in[0,1)$ we can write in compact form
\begin{equation}\label{eq:Ea}
\E_{\al}(s,q)=\frac{1}{\sum_{n\geq 1}\omega^{\al}_{n}q^{n-1}}\sum_{n\geq 2}\omega_{n}^{\al}(2-q)q^{n-2}A_{n}(s).
\end{equation}

\begin{proposition}\label{prop:E}
For any $\al\in[0,1)$ and $q\in(0,1)$ the function $\mathcal E_{\alpha}(q,s)$ is monotonically decreasing and convex in $s$ for any $s\in(0,1)$.
\end{proposition}
\proof
Since $A'_{2}(s)=0$, for  any  $\alpha\in[0,1)$
$$
\frac{d}{ds}\E_{\alpha}(s,q)= \frac{1}{\sum_{n\geq 1}\omega^{\al}_{n}q^{n-1}}\sum_{n\geq 3}\omega^{\al}_{n}q^{n-2}(2-q)A'_{n}(s)
$$
and
$$
\frac{d^{2}}{ds^{2}}\E_{\alpha}(s,q)= \frac{1}{\sum_{n\geq 1}\omega^{\al}_{n}q^{n-1}}\sum_{n\geq 4}\omega^{\al}_{n}q^{n-2}(2-q)A''_{n}(s).
$$
Since $\omega_{n}^{\al}>0$ and, for point $iii)$ of lemma \ref{lem:propA}, $A'_{n}<0$ and $A''_{n}>0$, the thesis follows. 

\qed

\begin{corollary}\label{cor:Epos}
For any $\alpha\in[0,1)$ the function $\E_{\alpha}(s,q)> 0$ for any $s\in(0,1)$, $q\in(0,1)$.
\end{corollary}
\proof
From the previous proposition
$$
\E_{\alpha}(s,q)\geq \E_{\alpha}(1,q)=\left\{\begin{array}{ll}
\frac{2-q}{q}\left(\frac{1-(1-q)^{\alpha}-\alpha q}{1-(1-q)^{\alpha}}\right), &0<\alpha<1\\
\frac{2-q}{q}\left(1+\frac{q}{\log(1-q)}\right), &\alpha=0.
\end{array}\right.
$$
Since $(1-q)^{\alpha}<1-\alpha q$ and $-\log(1-q)>q$, it follows $ \E_{\alpha}(1,q)>0$ for any $q\in(0,1)$.

\qed

Consider now the derivative of $\E_{\alpha}(s,q)$ with respect to the variable $q$. For $\alpha\in[0,1)$

\begin{equation}\label{eq:Dedq}
\begin{split}
\partial_{q}\E_{\alpha}(s,q)&=\frac{1}{\left(\sum_{n\geq 1}\omega^{\al}_{n}q^{n-1} \right)^{2}}\cdot\left\{\left(\sum_{n\geq 2}\omega^{\al}_{n}\Big((n-2)(2-q)-q\Big)q^{n-3}A_{n}(s)\right)\left(\sum_{n\geq 1}\omega^{\al}_{n}q^{n-1}\right)\right.\\
&\hspace{100pt}\left.-\left(\sum_{n\geq 2}\omega^{\al}_{n}(2-q)q^{n-2}A_{n}(s)\right)\left(\sum_{n\geq 1}\omega^{\al}_{n}(n-1)q^{n-2}\right)\right\}\\
&=\frac{1}{\left(\sum_{n\geq 1}\omega^{\al}_{n}q^{n-1} \right)^{2}}\cdot\sum_{n\geq 2,m\geq 1}\omega^{\al}_{n}\omega^{\al}_{m}\Big((2-q)(n-m)-2\Big)A_{n}(s)q^{(n-2)+(m-2)}.
\end{split}
\end{equation}
Henceforth, let us denote
$$
C_{\al}(q)\bydef\left(\sum_{n\geq 1}\omega^{\al}_{n}q^{n-1} \right)^{-2}.
$$
In the sequel we prove that the function $\partial_{q}\E_{\alpha}(s,q)$ is monotonic and convex in $s$. For that let us  introduce the quantity
$$
H_{\alpha}(n,m)\bydef \omega^{\al}_{n+2}\omega^{\al}_{m}(n-m)
$$
and we prove the following.

\begin{lemma}\label{lem:Hnm}
Let $n,m\geq 1$ and  $n>m$. Then, for any $\al\in[0,1)$,
$$
H_{\alpha}(n,m)+H_{\alpha}(m,n)>0.
$$
\end{lemma}
\proof
Case  $\al=0$. 
 $$
  H_{0}(n,m)+H_{0}(m,n)=\left(\frac{1}{n+2}\frac{1}{m}-\frac{1}{m+2}\frac{1}{n}\right)(n-m)=2\frac{(n-m)^{2}}{nm(n+2)(m+2)}>0.
  $$
 Let be now $\al\in(0,1)$.
For any integer $t$, we remind the relations
$$
{\al \choose t+1}=\frac{\al-t}{t+1}{\al \choose t},\quad {\al \choose t+2}=\frac{\al-(t+1)}{t+2}\frac{\al-t}{t+1}{\al \choose t}.
$$
Then
\begin{equation}
\begin{split}
H_{\al}(n,m)+H_{\al}(m,n)=&\left(\frac{\alpha-(n+1)}{n+2}\frac{\alpha-n}{n+1}-\frac{\alpha-(m+1)}{m+2}\frac{\alpha-m}{m+1}\right)\omega^{\al}_{n}\omega^{\al}_{m}(n-m)\\
=&\left(\frac{(n+1)-\alpha}{n+2}\frac{n-\alpha}{n+1}-\frac{(m+1)-\alpha}{m+2}\frac{m-\alpha}{m+1}\right)\omega^{\al}_{n}\omega^{\al}_{m}(n-m)\ .\\
\end{split}
\end{equation}
Since $\frac{K-\alpha}{K+1}>\frac{P-\alpha}{P+1}>0$ whenever $K>P\geq 1$, it follows $H_{\al}(n,m)+H_{\al}(m,n)>0$ for any $n>m$.
\qed

\begin{proposition}\label{prop:Eq}
For any $\alpha\in[0,1)$ and $q\in(0,1)$
\begin{itemize}
\item[i)] $$
\frac{d}{ds}\Big( \partial_{q}\E_{\alpha}(s,q)\Big)<0,\quad \forall  s\in(0,1)$$
\item[ii)]
$$\frac{d^{2}}{ds^{2}}\Big( \partial_{q}\E_{\alpha}(s,q)\Big)>0,\quad \forall s\in(0,1).$$
\end{itemize}
\end{proposition}
\proof

$i)$ From \eqref{eq:Dedq}, since $A'_{2}(s)=0$,
\begin{equation*}
\begin{split}
\frac{d}{ds}\left( \partial_{q}\E_{\alpha}(s,q)\right)&=C_{\al}(q)\sum_{\substack{n\geq 3, m\geq 1}}\omega^{\al}_{n}\omega^{\al}_{m}\Big((2-q)(n-m)-2\Big)A'_{n}(s)q^{(n-2)+(m-2)}\\
&=C_{\al}(q)\sum_{p\geq 4}\left[\sum_{\substack{n+m=p\\n\geq 3, m\geq 1}}\omega^{\al}_{n}\omega^{\al}_{m}\Big((2-q)(n-m)-2\Big)A'_{n}(s)\right]q^{p-4}.
\end{split}
\end{equation*}
Denote by 
$$
Coef(s,q,\alpha,p)\bydef\sum_{\substack{n+m=p\\n\geq 3, m\geq 1}}\omega^{\al}_{n}\omega^{\al}_{m}\Big((2-q)(n-m)-2\Big)A'_{n}(s).
$$
We prove that $Coef(s,q,\alpha,p)<0$ for any $s\in(0,1)$, $q\in(0,1)$ and $\alpha\in[0,1)$ and $p\geq 4$. Compute first the derivative in $q$:
\begin{equation*}
\begin{split}
\frac{d}{dq}Coef(s,q,\alpha,p)=&-\sum_{\substack{n+m=p\\n\geq 3, m\geq 1}}\omega^{\al}_{n}\omega^{\al}_{m}(n-m)A'_{n}(s).\\
\end{split}
\end{equation*}
For point $iii)$ in lemma \ref{lem:propA} $A'_{n}(s)<0$, then if  $p\leq 6$ the above sum is obviously positive, being $n\geq m$ in any contribution. For $p\geq7$, we write
$$
\frac{d}{dq}Coef(s,q,\alpha,p)=-\omega^{\al}_{p-1}\omega^{\al}_{1}(p-2)A'_{p-1}(s)-\omega^{\al}_{p-2}\omega^{\al}_{2}(p-4)A'_{p-2}(s)-\sum_{\substack{n+m=p\\n\geq 3, m\geq 3}}\omega^{\al}_{n}\omega^{\al}_{m}(n-m)A'_{n}(s).\\
$$
The first and the second terms in the right hand side are positive then, collecting the contributions to the sum due to the couples $(n,m)$ and $(m,n)$,   and for  point $iv)$ in lemma \ref{lem:propA}, it holds
$$
\frac{d}{dq}Coef(s,q,\alpha,p)>-\sum_{\substack{n+m=p\\n>m\geq 3}}\omega^{\al}_{n}\omega^{\al}_{m}(n-m)(A'_{n}(s)-A'_{m}(s))>0.
$$
Therefore,  $\frac{d}{dq}Coef(s,q,\alpha,p)>0$ for any $p\geq 4$ and $s\in(0,1)$, $q\in(0,1)$,  $\al\in[0,1)$. It means that 
$$
Coef(s,q,\alpha,p)<Coef(s,1,\alpha,p).
$$
Now,
\begin{equation*}
\begin{split}
Coef(s,1,\alpha,p)&
=\sum_{\substack{n+m=p\\n\geq 3, m\geq 1}}\omega^{\al}_{n}\omega^{\al}_{m}\Big(n-m-2\Big)A'_{n}(s)=\sum_{\substack{n+m=p-2\\n\geq 1, m\geq 1}}\omega^{\al}_{n+2}\omega^{\al}_{m}\Big(n-m\Big)A'_{n+2}(s)\\
&=\sum_{\substack{n+m=p-2\\n>m\geq 1}}\left[\omega^{\al}_{n+2}\omega^{\al}_{m}\Big(n-m\Big)A'_{n+2}(s)+\omega^{\al}_{m+2}\omega^{\al}_{n}\Big(m-n\Big)A'_{m+2}(s)\right].
\end{split}
\end{equation*}
Since $\omega^{\al}_{n+2}\omega^{\al}_{m}>0$ and, again for property $iii)$ of lemma \ref{lem:propA}, $A'_{n+2}(s)<A'_{m+2}(s)$, it follows

\begin{equation*}
\begin{split}
Coef(s,1,\al,p)&<\sum_{\substack{n+m=p\\n>m\geq 1}}\left[\omega^{\al}_{n+2}\omega^{\al}_{m}\Big(n-m\Big)+\omega^{\al}_{m+2}\omega^{\al}_{n}\Big(m-n\Big)\right]A'_{m+2}(s)\\
&=\sum_{\substack{n+m=p\\n>m\geq 1}}(H_{\al}(n,m)+H_{\al}(m,n))A'_{m+2}(s)<0
\end{split}
\end{equation*}
in light of lemma \ref{lem:Hnm}.
Therefore all the coefficients $Coef(s,q,\al,p)<0$ and  $\frac{d}{ds}\left( \partial_{q}\E_{\alpha}(s,q)\right)<0$.

$ii)$
Consider now the second derivative
\begin{equation*}
\begin{split}
&\frac{d^{2}}{ds^{2}}\Big(\partial_{q}\E(s,q)\Big)=C_{\al}(q)
\sum_{p\geq 5}\sum_{\substack{n+m=p\\n\geq 4, m\geq 1}}\omega^{\al}_{n}\omega^{\al}_{m}\Big((2-q)(n-m)-2\Big)A''_{n}(s)q^{p-4}\ .
\end{split}
\end{equation*}
Arguing as before, we realise that $q=1$ represents the worst case and, since $A''_{n}(s)>0$,  it follows 
$$
\sum_{\substack{n+m=p\\n\geq 4, m\geq 1}}\omega^{\al}_{n}\omega^{\al}_{m}\Big((n-m)-2\Big)A''_{n}(s)>0,\quad \forall s\in(0,1), q\in(0,1),\alpha\in(0,1), p\geq 5.
$$

\qed

Up to now we investigated the monotonicity and the convexity of $\partial_{q}\E_{\al}(q,s)$ without mentioning about the value of $\partial_{q}\E_{\al}(q,s)$. To the purpose of proving the monotonicity of the apsidal angle, it would be enough that $\partial_{q}\E_{\al}(q,s)$ is positive for any $s$. Unfortunately that is not the case, indeed  for any $\alpha$ there is a region of $q$ where  $\partial_{q}\E_{\al}(q,s)$ is not everywhere positive. Nevertheless,  we prove now that  $\partial_{q}\E_{\al}(q,s)>0$  at $s=1/2$ for any $\alpha$ and $q$. For the previous result it follows that $\partial_{q}\E_{\al}(q,s)>0$ for any $s\in(0,1/2)$, for any $q$ and $\al$.  

\begin{proposition}\label{pro:eq12}
For any $q\in(0,1)$ and $\alpha\in[0,1)$
$$
\partial_{q}\E_{\alpha}\left(\frac{1}{2},q\right)>0.
$$
\end{proposition}

\proof

\begin{equation*}
\begin{split}
\partial_{q}\E_{\al}\left(\frac{1}{2},q\right)&=C_{\al}(q)\sum_{p\geq 3}\sum_{\substack{n+m=p\\n\geq 2, m\geq 1}}\omega^{\al}_{n}\omega^{\al}_{m}\Big((2-q)(n-m)-2\Big)A_{n}(1/2)q^{p-4}\\
&=C_{\al}(q)\sum_{p\geq 3}\left[ \sum_{\substack{n+m=p\\n\geq 2, m\geq 1}}\omega^{\al}_{n}\omega^{\al}_{m}\Big(2(n-m-1)\Big)A_{n}(1/2)q^{p-4}\right.\\
&\hspace{50pt}\left.-\sum_{\substack{n+m=p\\n\geq 2, m\geq 1}}\omega^{\al}_{n}\omega^{\al}_{m}\Big(n-m\Big)A_{n}(1/2)q^{p-3}\right].\\
\end{split}
\end{equation*}
When $p=3$ the first sum does not contribute, then we can write
\begin{equation*}
\begin{split}
\partial_{q}\E\left(\frac{1}{2},q\right)&=C_{\al}(q)\sum_{p\geq 4}\left[ \sum_{\substack{n+m=p\\n\geq 2, m\geq 1}}\omega^{\al}_{n}\omega^{\al}_{m}2(n-m-1)A_{n}(1/2)-\sum_{\substack{n+m=p-1\\n\geq 2, m\geq 1}}\omega^{\al}_{n}\omega^{\al}_{m}\Big(n-m\Big)A_{n}(1/2)\right]q^{p-4}.
\end{split}
\end{equation*}
Denote by 
$$
Q_{p}(q)\bydef\sum_{\substack{n+m=p\\n\geq 2, m\geq 1}}\omega^{\al}_{n}\omega^{\al}_{m}2(n-m-1)A_{n}(1/2)-\sum_{\substack{n+m=p-1\\n\geq 2, m\geq 1}}\omega^{\al}_{n}\omega^{\al}_{m}(n-m)A_{n}(1/2)
$$
 the coefficient of $q^{p-4}$.
For $p\geq 4$
 \begin{equation}
\begin{split}
Q_{p}(q)&=\sum_{\substack{n+m=p-1\\n\geq 1, m\geq 1}}\omega^{\al}_{n+1}\omega^{\al}_{m}2(n-m)A_{n+1}(1/2)-\sum_{\substack{n+m=p-1\\n\geq 2, m\geq 1}}\omega^{\al}_{n}\omega^{\al}_{m}(n-m)A_{n}(1/2)\\
&=\omega^{\al}_{2}\omega^{\al}_{p-2}2(3-p)A_{2}(1/2)+\sum_{\substack{n+m=p-1\\n\geq 2, m\geq 1}}\omega^{\al}_{n}\omega^{\al}_{m}(n-m)\Big(2\frac{\omega^{\al}_{n+1}}{\omega^{\al}_{n}}A_{n+1}(1/2)-A_{n}(1/2)\Big)\\
&=\omega^{\al}_{2}\omega^{\al}_{p-2}2(3-p)A_{2}(1/2)+\omega^{\al}_{p-2}\omega^{\al}_{1}(p-3)\left(2\frac{\omega^{\al}_{p-1}}{\omega^{\al}_{p-2}}A_{p-1}(1/2)-A_{p-2}(1/2) \right)\\
&\quad+\sum_{\substack{n+m=p-1\\n\geq 2, m\geq 2}}\omega^{\al}_{n}\omega^{\al}_{m}(n-m)\left(\frac{2\omega^{\al}_{n+1}}{\omega^{\al}_{n}}A_{n+1}(1/2)-A_{n}(1/2)\right).\\
\end{split}
\end{equation}
Reminding the definition \eqref{eq:Kn} of $K_{n}^{\al}(s)$ and collecting the contributions due to the couple $(n,m)$ and $(m,n)$, it results
 \begin{equation}
\begin{split}
Q_{p}(q)&=\omega^{\al}_{p-2}(p-3)\left[2\omega^{\al}_{1}\frac{\omega^{\al}_{p-1}}{\omega^{\al}_{p-2}}A_{p-1}(1/2)-\omega^{\al}_{1}A_{p-2}(1/2)-2\omega^{\al}_{2}A_{2}(1/2) \right]\\
&\quad +\sum_{\substack{n+m=p-1\\n>m\geq 2}}\omega^{\al}_{n}\omega^{\al}_{m}(n-m)\Big( K^{\al}_{n}(1/2)-K^{\al}_{m}(1/2)\Big).
\end{split}
\end{equation}
For point $ii)$ in lemma \ref{lem:propKn}, the last sum is positive.
Moreover, for $\al\in(0,1)$, the function
 \begin{equation}
\begin{split}
\eta^{\al}(p)&\bydef2\omega^{\al}_{1}\frac{\omega^{\al}_{p-1}}{\omega^{\al}_{p-2}}A_{p-1}(1/2)-\omega^{\al}_{1}A_{p-2}(1/2)-2\omega^{\al}_{2}A_{2}(1/2) \\
&=\frac{2\al(p-2-\al)}{p-1}A_{p-1}(1/2)-\al A_{p-2}(1/2) -\al(1-\al)
\end{split}
\end{equation}
is such that
$$
\eta^{\al}(4)=0,\quad  \dot\eta^{\al}(p)=\frac{2\al(2+2\al)}{(p-1)^{2}2^{p-2}}\Big(2^{p-2}-1-(p-1)\log(2)\Big)>0 \quad \forall p\geq 4.
$$
A similar results holds for $\eta^{0}(p)$.
Thus $Q_{p}(q)>0$ for any $p\geq 4$ and $\alpha\in[0,1)$ and the thesis follows.

\qed

\section{Limits of $ \Delta_\alpha\theta(u)$ as $\ell\to 0$ and $\ell\to \ell_{max}$}\label{sec:limits}

This section concerns the limits of the apsidal angle for near-circular and near-radial orbits, that is when $\ell\to \ell_{max}$ and $\ell\to 0$.

Although these limit values are well know in the literature, we present the result for the sake of completeness and because it descends directly from the above computations.
\begin{theorem}
For any $\alpha\in[0,1)$
$$
\lim_{\ell\to 0}\Delta_{\alpha}\theta(u)=\frac{\pi}{2-\alpha}, \quad\quad  \lim_{\ell\to \ell_{\max}}\Delta_{\alpha}\theta(u)=\frac{\pi}{\sqrt{2-\alpha}}.
$$
\end{theorem}  
\proof
According to the discussion in section \ref{sec:aps-ang}, the limit values for the apsidal angle are given as the limit for $q\to 0$ and $q\to 1$ of the integral \eqref{eq:deltatheta}. In particular, $q\to 0$ corresponds to $\ell\to \ell_{max}$ and $q\to 1$ to $\ell\to 0$. Since $\E_{\alpha}(s,q)>0$ for any $\alpha\in[0,1), q\in(0,1) , s\in(0,1)$, see Corollary \ref{cor:Epos}, for the Lebesgue's dominated convergence theorem,
$$
\lim_{q\to 0}\Delta_{\alpha}\theta=\int_{0}^{1}\frac{1}{\sqrt{s(1-s)}}\frac{1}{\sqrt{1+\lim_{q\to 0}\E_{\alpha}(s,q)}}\, ds=\int_{0}^{1}\frac{1}{\sqrt{s(1-s)}}\frac{1}{\sqrt{1+(1-\alpha)}}\, ds=\frac{\pi}{\sqrt{2-\alpha}}
$$
and
$$
\lim_{q\to 1}\Delta_{\alpha}\theta=\int_{0}^{1}\frac{1}{\sqrt{s(1-s)}}\frac{1}{\sqrt{1+\lim_{q\to 1}\E_{\alpha}(s,q)}}\, ds=\int_{0}^{1}\frac{1}{\sqrt{s(1-s)}}\frac{1}{\sqrt{1+\frac{1}{s}\frac{s^{\al}-s}{1-s}}}\, ds=\frac{\pi}{2-\alpha}.
$$
\qed

\section{The derivative of  $ \Delta_\alpha\theta(u)$}\label{sec:derivative}

We now examine  the derivative of the apsidal angle as function of the angular momentum $\ell$. That is equivalent to compute the derivative of $\Delta_{\alpha}\theta(u)$ given in \eqref{eq:deltatheta} with respect to $q$. From  proposition \ref{prop:Eq} and   \ref{pro:eq12} 
$$
|\partial_{q}\E_{\al}(s,q)|<\partial_{q}\E_{\al}(0,q)=M_{\alpha}(q)
$$
where 
$$
M_{\alpha}(q)=\left\{\begin{array}{ll}
{\begin{array}{l}
a^2\,\left(1-q\right)^{a}\,q^3+\left(\left(2-a\right)\,\left(1-q
 \right)^{2\,a}+\left(-2\,a^2+a-4\right)\,\left(1-q\right)^{a}+2
 \right)\,q^2\\
 +\left(-4\,\left(1-q\right)^{2\,a}+8\,\left(1-q\right)^{
 a}-4\right)\,q+2((1-q)^{\al}-1)^{2} 
 \end{array}}\over{{\displaystyle \left(\left(1-q\right)^{a}-1\right)^2\,\left(q-1\right)^2\,q^
 2}}, &\al\in(0,1)\\
 \\
{\displaystyle \frac {{q}^{3}-\ln  \left( 1-q \right) {q}^{2}-2\,{q}^{2}+2\, \left( 
\ln  \left( 1-q \right)  \right) ^{2}{q}^{2}-4\, \left( \ln  \left( 1-
q \right)  \right) ^{2}q+2\, \left( \ln  \left( 1-q \right)  \right) ^
{2}}{ \left( \ln  \left( 1-q \right)  \right) ^{2}{q}^{2} \left( 1-2\,
q+{q}^{2} \right)} }
& \al=0
\end{array}\right.
$$
is $\mathcal C^{1}((0,1))$ for any $\alpha$. Therefore, for any $q\in(0,1)$ we can pass the derivative under the integral sign and write
\begin{equation}\label{eq:der}
\frac{d}{dq}\Delta_{\al}\theta(u)=\int_{0}^{1}\frac{1}{\sqrt{s(1-s)}}\partial_{q}\left( \frac{1}{\sqrt{1+\E_{\al}(s,q)}}\right)\,ds=-\frac{1}{2}\int_{0}^{1}\frac{1}{\sqrt{s(1-s)}} \frac{\partial_{q}\E_{\al}(s,q)}{\left(1+\E_{\al}(s,q)\right)^{\frac{3}{2}}}\,ds.
\end{equation}
The following analysis aims at proving that the integral in \eqref{eq:der} is positive for any $q\in(0,1)$ and $\alpha\in[0,1)$. It would mean that the apsidal angle is monotonically increasing as function of the angular momentum $\ell$ for any $\alpha\in[0,1)$.

Denote by
$$
I_{\al}(q)\bydef\int_{0}^{1}\frac{1}{\sqrt{s(1-s)}} \frac{\partial_{q}\E_{\al}(s,q)}{\left(1+\E_{\al}(s,q)\right)^{\frac{3}{2}}}\,ds.
$$
Taking advantage from the symmetry in $s$, it holds
$$
I_{\al}(q)=\int_{0}^{1/2}\frac{1}{\sqrt{s(1-s)}}\left( \frac{\partial_{q}\E_{\al}(s,q)}{\left(1+\E_{\al}(s,q)\right)^{\frac{3}{2}}}+ \frac{\partial_{q}\E_{\al}(1-s,q)}{\left(1+\E_{\al}(1-s,q)\right)^{\frac{3}{2}}}\right)\, ds.
$$
For propositions \ref{prop:Eq}, \ref{pro:eq12}, the first term of the integrand is positive, but, as said earlier, $\partial_{q}\E_{\al}(1-s,q)$ could be negative.
A sufficient condition for $I_{\al}(q)>0$ is that 
\begin{equation}\label{eq:disug}
\frac{\partial_{q}\E_{\al}(s,q)}{\left(1+\E_{\al}(s,q)\right)^{\frac{3}{2}}}>-\frac{\partial_{q}\E_{\al}(1-s,q)}{\left(1+\E_{\al}(1-s,q)\right)^{\frac{3}{2}}},\quad \forall s\in(0,1/2).
\end{equation}
 The convexity of $\partial_{q}\E_{\al}(s,q)$ in the variable $s$ implies that $\partial_{q}\E_{\al}(s,q)> \left| \partial_{q}\E_{\al}(1-s,q)\right|$ for any $s\in(0,1/2)$. On the other side, from proposition \ref{prop:E},  $\E_{\al}(s,q)>\E_{\al}(1-s,q)$ hence the validity of  \eqref{eq:disug} can not be achieved from the global behaviour of $\E_{\al}$ and $\partial_{q} \E_{\al}$, rather  it deserves a carefully analysis.   

Consider the function
$$
\mathcal I_{\al}(s,q)\bydef\partial_{q}\E_{\al}(s,q)\left(1+\E_{\al}(1-s,q)\right)^{\frac{3}{2}}+\partial_{q}\E_{\al}(1-s,q)\left(1+\E_{\al}(s,q)\right)^{\frac{3}{2}}.
$$
Clearly, the statement $\mathcal I_{\al}(s,q)>0$ for any $s\in(0,1/2)$ is equivalent to \eqref{eq:disug}, then we aim at proving that $\mathcal I_{\al}(s,q)$ is positive. The next lemma allows to replace the power $\frac{3}{2}$ by $2$.
\begin{lemma}\label{lem:3/2}
Let $A\geq C>0$, $B,D>0$.
If 
$$
AB^{2}-CD^{2}>0\Rightarrow AB^{\frac{3}{2}}-CD^{\frac{3}{2}}>0.
$$
\end{lemma}
\proof
$$
AB^{2}>CD^{2}\Rightarrow A^{\frac{3}{4}}B^{\frac{3}{2}}>C^{\frac{3}{4}}D^{\frac{3}{2}}\Rightarrow AB^{\frac{3}{2}}>C^{\frac{3}{4}}A^{\frac{1}{4}}D^{\frac{3}{2}}>CD^{\frac{3}{2}}.
$$
\qed

For those $s$ where $\partial_{q}\E_{\al}(q,1-s)\geq 0$ it holds $\mathcal I_{\al}(s,q)>0$ and \eqref{eq:disug} is true. While for those $s$ where $\partial_{q}\E_{\al}(q,1-s)<0$ it holds $\partial_{q}\E_{\al}(q,s)>-\partial_{q}\E_{\al}(q,1-s)$, then for the lemma, we can replace the above function $\mathcal I_{\al}(s,q)$ by 
\begin{equation}\label{eq:Ial}
\mathcal I_{\al}(s,q)\bydef\partial_{q}\E_{\al}(q,s)\left(1+\E_{\al}(q,1-s)\right)^{2}+\partial_{q}\E_{\al}(q,1-s)\left(1+\E_{\al}(q,s)\right)^{2}
\end{equation}
and we investigate whether the new defined $\mathcal I_{\al}(s,q)$ is positive.

From \eqref{eq:Dedq}  it descends, for $\alpha\in(0,1)$,
\begin{equation}\label{eq:Ia2}
\begin{split}
&\mathcal I_{\al}(q,s)=C_{\al}(q)\sum_{\substack{n\geq 2, m\geq 1}}\omega^{\al}_{n}\omega^{\al}_{m}\Big((2-q)(n-m)-2\Big)\left[\begin{array}{l}
A_{n}(s)\left(1+\E_{\al}(1-s,q)\right)^{2}\\
+A_{n}(1-s)\left(1+\E_{\al}(s,q)\right)^{2}\end{array}\right]q^{(n-2)+(m-2)}.
\end{split}
\end{equation}
Algebraic manipulations lead to the expansion 
\begin{equation}\label{eq:mI}
\mathcal I_{\al}(s,q)=C_{\al}(q)\sum_{p\geq 4}Q_{\al}(p,s)q^{p-4}
\end{equation}
where 
 \begin{equation}\label{eq:Qps1}
\begin{split}
Q_{\al}(p,s)=&\sum_{\substack{n\geq 2, m\geq 1\\n+m=p}}\omega^{\al}_{n}\omega^{\al}_{m}2(n-m-1)\left[A_{n}(s)\Big(1+\E_{\al}(1-s,q)\Big)^{2}+A_{n}(1-s)\Big(1+\E_{\al}(s,q)\Big)^{2} \right]\\
&-\sum_{\substack{n\geq 2, m\geq 1\\n+m=p-1}}\omega^{\al}_{n}\omega^{\al}_{m}(n-m)\left[A_{n}(s)\Big(1+\E_{\al}(1-s,q)\Big)^{2}+A_{n}(1-s)\Big(1+\E_{\al}(s,q)\Big)^{2} \right].
\end{split}
\end{equation}
Note that the  quantities in square brackets are positive and increasing in $n$ and it can be proved that both the sums appearing in $Q_{\al}(p,s)$ are positive for any $p\geq 4$. However,   $Q_{\al}(p,s)$ has not a definite sign as $q$ ranges in $(0,1)$ then we can not conclude about the sign of $\mathcal I_{\al}(s,q)$ looking at the sign of each of the  $Q_{\al}(p,s)$'s. Therefore we have to  extrapolate the dependence on $q$ from each of the $\E_{\al}(s,q)$ and $\E_{\al}(1-s,q)$.
By substituting  the expansion \eqref{eq:Ea} in \eqref{eq:Ia2}, 
straightforward computations give
\begin{equation}\label{eq:Iaserie}
\begin{split}
\mathcal I_{\al}(q,s)=C_{\al}^{2}(q)\sum_{\substack{n\geq 2, m\geq 1\\k_{1},k_{2}\geq 2}}&\omega^{\al}_{n}\omega^{\al}_{m}\omega^{\al}_{k_{1}}\omega^{\al}_{k_{2}}\Big((2-q)(n-m)-2\Big)q^{(n-2)+(m-2)+(k_{1}-2)+(k_{2}-2)}\cdot\\
&\cdot\left[\begin{array}{l}
A_{n}(s)\left(\frac{\omega^{\al}_{k_{1}-1}}{\omega^{\al}_{k_{1}}}+(2-q)A_{k_{1}}(1-s)\right)\left(\frac{\omega^{\al}_{k_{2}-1}}{\omega^{\al}_{k_{2}}}+(2-q)A_{k_{2}}(1-s)\right)\\
+A_{n}(1-s)\left(\frac{\omega^{\al}_{k_{1}-1}}{\omega^{\al}_{k_{1}}}+(2-q)A_{k_{1}}(s)\right)\left(\frac{\omega^{\al}_{k_{2}-1}}{\omega^{\al}_{k_{2}}}+(2-q)A_{k_{2}}(s)\right)
\end{array}\right].
\end{split}
\end{equation}
Note that in the series there is no the term $q^{-1}$. Indeed the only possibility for this term to appear is for the choice $(n,m,k_{1},k_{2})=(2,1,2,2)$ but the first bracket results $-q$. Also the coefficient of $q^{0}$ is equal to zero. 
Therefore, collecting on the powers of $q$, we can write 
$$
\mathcal I_{\al}(q,s)=C^{2}_{\al}(q)\sum_{p\geq 1}Coef_{p}(\al,s)q^{p}
$$
where, for any $p$, only a finite number of quart-uple  $(n,m,k_{1},k_{2})$ contributes to  $Coef_{p}(\al,s)$.
Here we list the first few of these coefficients, for the case $\alpha\in(0,1)$.

$$
Coef_{1}(\al,s)=\frac{1}{9} \left( 2-\al \right) ^{2} \left( (5+\al){s}^{2}\al-(5+\al)s+\al+2 \right) 
$$
$$
Coef_{2}(\al,s)=\frac{1}{18}\al^{4}(1-\al)(7-4\al)(2-\al)^{3} \left( (5+\al){s}^{2}\al-(5+\al)s+\al+2 \right) 
$$
$$
Coef_{3}(\al,s)=\frac{1}{90}a^{4}(2-a)^{3}(1-a)\left(\begin{array}{l} \left( {a}^{3}-5\,{a}^{2}-17\,a+69 \right) {s}^{4}+
\left( 10\,{a}^{2
}-138-2\,{a}^{3}+34\,a \right) {s}^{3}+\\
 \left( 482-346\,a+15\,{a}^{2}+
23\,{a}^{3} \right) {s}^{2}+ \\
\left( -22\,{a}^{3}-20\,{a}^{2}-413+329\,
a \right) s+\\
21\,{a}^{3}-84\,a-39\,{a}^{2}+156\end{array}\right)
$$
$$
Coef_{4}(\al,s)=\frac{1}{540}\al^{4}(1-\al)(2-\al)^{3}(9-4a)\left(\begin{array}{l}
 \left( 3\,{a}^{3}-15\,{a}^{2}-51\,a+207 \right) {s}^{4}+\\
  \left( -6\,{
a}^{3}+30\,{a}^{2}+102\,a-414 \right) {s}^{3}+\\
 \left( 29\,{a}^{3}-5\,{
a}^{2}-428\,a+746 \right) {s}^{2}+ \\
\left( -26\,{a}^{3}-10\,{a}^{2}+377
\,a-539 \right) s+\\
23\,{a}^{3}-47\,{a}^{2}-92\,a+188
\end{array}
\right).
$$
The first two coefficients are clearly positive for any $s\in(0,1/2)$ and $\al\in(0,1)$ and the same it holds for $Coef_{3}(\al,s)$, $Coef_{4}(\al,s)$. 
By symbolic computation and numerical visualisation it appears evident that even for larger $p$ the coefficients $Coef_{p}(\al,s)$ are positive and, therefore, that $\mathcal I_{\al}(q,s)>0$. On the other side, although rigorous,  the check of the positivity of any large but  finite number of coefficients does not provide a proof for $\mathcal I_{\al}(q,s)>0$, unless a proof that $Coef_{p}(\al,s)>0$ for any $p>p^{\star}$ is provided.

That is exactly what we are presenting in the next section, where the monotonicity of the apsidal angle is proved in the case $\al=0$.

\section{The monotonicity of the apsidal angle in the logarithmic potential case}\label{sec:proof}

We aim at proving that, for any fixed value of the energy,  the apsidal angle $\Delta_{0}\theta(u)$ is monotonically increasing as function of the angular momentum $\ell$. According to the discussion of the previous sections, it is enough to prove that 
the function
\begin{equation}\label{eq:Ia0}
\mathcal I_{0}(s,q)=\partial_{q}\E_{0}(s,q)\left(1+\E_{0}(1-s,q)\right)^{2}+\partial_{q}\E_{0}(1-s,q)\left(1+\E_{0}(s,q)\right)^{2}
\end{equation}
is positive for any $s\in(0,1/2)$ and $q\in(0,1)$.

The proof is split in different steps and analytical arguments  are sometimes  combined with rigorous numerical computations based on interval arithmetics, \cite{MR2652784}.  All the numerical computation have been  performed in {\it Matlab} with the interval arithmetic package {\it Intlab}. \cite{Ru99a}.  That assures that the  results we obtain are reliable in the strict mathematical sense. 

First we show that $\partial_{q}\E_{0}(1-s,q)>0$ for any $s\in(0,1/2)$ for any $q\in[0.9,1)$, yielding $\mathcal I_{0}(s,q)>0$ for any $s\in(0,1/2)$, $q\in[0.9,1)$. 

\begin{proposition}\label{prop:0.9}
For any $q\in[0.9,1)$, $s\in(0,1/2)$, $\partial_{q}\E_{0}(1-s,q)>0$. 
\end{proposition}
\proof
From proposition \ref{prop:Eq}, it is enough to show that $\lim_{s\to 1}\partial_{q}\E_{0}(s,q)>0$ for $q\in[0.9,1)$.

It holds
$$
F(q)\bydef\lim_{s\to 1}\partial_{q}\E_{0}(s,q)=-\frac{2}{q^{2}}-\frac{1}{\log(1-q)}+\frac{2-q}{(1-q)\log^{2}(1-q)}.
$$
Compute
$$
F'(q)=\frac{4(1-q)^{2}\log^{3}(1-q)+q^{4}\log(1-q)+2q^{3}(2-q)}{q^{3}(1-q)^{2}\log^{3}(1-q)}=:\frac{N(q)}{D(q)}.
$$
The denominator $D(q)$ is negative for any $q\in(0,1)$, while the numerator 
$$
N(q)<-4(1-q)^{2}q^{3}+q^{4}\log(1-q)+2q^{3}(2-q)=q^{4}\Big( \log(1-q)-4q+6\Big)=:q^{4}B(q).
$$
The function  $B(q)$ is decreasing in $q$ and, by  rigorous computation it holds $B(0.91)<-0.04$. Moreover, for any $q\in[0.9,0.91]$ it holds $N(q)<-0.2$.  Therefore $N(q)<0$  and $F(q)$ is increasing  for any $q\in[0.9,1)$. Since $F(0.9)\in[0.0394, 0.040]$ it follows that $F(q)>0$ for any $q\in[0.9,1)$.

\qed

From now on, we restrict our analysis to the case $q\in(0,0.9]$. For these values of $q$ the function $\partial_{q}\E_{0}(1-s,q)$ is not everywhere positive for $s\in(0,1/2)$ then the argument above adopted is not anymore valid to prove that $\mathcal I_{0}(s,q)>0$.

From \eqref{eq:mI}, we can separate
\begin{equation}\label{eq:separation}
\mathcal I_{0}(s,q)=C_{0}(q)\left(\mathcal I_{0}^{\boldsymbol f}(s,q)+\mathcal I_{0}^{\boldsymbol \infty}(s,q)\right)=C_{0}(q)\left(\sum_{p=4}^{10}Q_{0}(s,p)q^{p-4}+\sum_{p\geq 11}Q_{0}(s,p)q^{p}\right)
\end{equation}
and we prove separately, in section \ref{sec:infinite} and \ref{sec:finite}, that  the series $\mathcal I_{0}^{\boldsymbol \infty}(s,q)$  and the finite sum $\mathcal I_{0}^{\boldsymbol f}(s,q)$ are positive for any $s\in(0,1/2)$ and $q\in(0,0.9]$.  Reminding the definition \eqref{eq:Qps1} and \eqref{eq:w},  for the logarithmic potential case, we have
 \begin{equation}\label{eq:Qps0}
\begin{split}
Q_{0}(p,s)=&\sum_{\substack{n\geq 2, m\geq 1\\n+m=p}}\frac{1}{n}\frac{1}{m}2(n-m-1)\left(A_{n}(s)\Big(1+\E_0(1-s,q)\Big)^{2}+A_{n}(1-s)\Big(1+\E_0(s,q)\Big)^{2} \right)\\
&-\sum_{\substack{n\geq 2, m\geq 1\\n+m=p-1}}\frac{1}{n}\frac{1}{m}(n-m)\left(A_{n}(s)\Big(1+\E_0(1-s,q)\Big)^{2}+A_{n}(1-s)\Big(1+\E_0(s,q)\Big)^{2} \right).
\end{split}
\end{equation}
It is convenient to express the coefficients $Q_{0}(p,s)$ in a different way:
by a change of index $n\to n+1$ in the first sum, it follows
\begin{equation*}
\begin{split}
Q_{0}(p,s)=&\sum_{\substack{n\geq 1, m\geq 1\\n+m=p-1}}\frac{1}{n+1}\frac{1}{m}2(n-m)\left(A_{n+1}(s)\Big(1+\E_0(1-s,q)\Big)^{2}+A_{n+1}(1-s)\Big(1+\E_0(s,q)\Big)^{2} \right)\\
&-\sum_{\substack{n\geq 2, m\geq 1\\n+m=p-1}}\frac{1}{n}\frac{1}{m}(n-m)\left(A_{n}(s)\Big(1+\E_0(1-s,q)\Big)^{2}+A_{n}(1-s)\Big(1+\E_0(s,q)\Big)^{2} \right)\\
=&\frac{1}{2}\frac{1}{p-2}2(1-(p-2))\left( A_{2}(s)\Big(1+\E_0(1-s,q)\Big)^{2}+A_{2}(1-s)\Big(1+\E_0(s,q)\Big)^{2} \right)\\
&+\sum_{\substack{n\geq 2, m\geq 1\\n+m=p-1}}\frac{1}{n+1}\frac{1}{m}2(n-m)\left(A_{n+1}(s)\Big(1+\E_0(1-s,q)\Big)^{2}+A_{n+1}(1-s)\Big(1+\E_0(s,q)\Big)^{2} \right)\\
&-\sum_{\substack{n\geq 2, m\geq 1\\n+m=p-1}}\frac{1}{n}\frac{1}{m}(n-m)\left(A_{n}(s)\Big(1+\E_0(1-s,q)\Big)^{2}+A_{n}(1-s)\Big(1+\E_0(s,q)\Big)^{2} \right).
\end{split}
\end{equation*}
Collecting all the contributions and reminding the definition \eqref{eq:Kn}
it results
\begin{equation}\label{eq:Qps}
\begin{split}
Q_{0}(p,s)=&\sum_{\substack{n\geq 2, m\geq 1\\n+m=p-1}}\frac{1}{n}\frac{1}{m}(n-m)\left[
K_{n}(s)\Big(1+\E_0(1-s,q)\Big)^{2}+K_{n}(1-s)\Big(1+\E_0(s,q)\Big)^{2} \right]\\
&-\frac{p-3}{p-2}\Bigg( \Big(1+\E_0(1-s,q)\Big)^{2} +\Big(1+\E_0(s,q)\Big)^{2} \Bigg).
\end{split}
\end{equation}

\subsection{Analytical estimate of the tail elements }\label{sec:infinite}
The goal of this section is to prove that $Q_{0}(s,p)\geq 0$ for any $s\in(0,1/2)$, $q\in(0,1)$ and $p\geq 11$. That implies that $\mathcal I_{0}^{\boldsymbol \infty}(s,q)\geq 0$  for any $s\in(0,1/2)$ and $q\in(0,1)$.

In practice, we are going to prove that 
$$
\sum_{\substack{n\geq 2, m\geq 1\\n+m=p-1}}\frac{1}{n}\frac{1}{m}(n-m)K^{0}_{n}(s)>\frac{p-3}{p-2},\quad  \forall p\geq 11, s\in(0,1).
$$
A first estimates is the following:
\begin{lemma}\label{lem:bound1}
For any $p\geq 10$
$$
\sum_{\substack{n\geq 2, m\geq 1\\n+m=p-1}}\frac{1}{n}\frac{1}{m}(n-m)K^{0}_{n}(s)\geq\sum_{\substack{n\geq 2, m\geq 1\\n+m=p-1}}\frac{1}{n}\frac{1}{m}(n-m)\frac{n-1}{n+1},\quad \forall s\in(0,1).
$$
\end{lemma}
\proof
Set
$$
g(s)\bydef \sum_{\substack{n\geq 2, m\geq 1\\n+m=p-1}}\frac{1}{n}\frac{1}{m}(n-m)K^{0}_{n}(s).
$$
Assume $p$ is odd, $p=2l+1$.
Then
\begin{equation*}
\begin{split}
g(s)&=\sum_{\substack{n=l+1\\(m=2l-n)}}^{2l-2}\frac{1}{n}\frac{1}{m}(n-m)\Big(K^{0}_{n}(s)-K^{0}_{m}(s)\Big)+\frac{2l-2}{2l-1}K^{0}_{2l-1}(s)\\
&=\sum_{\substack{n=l+1\\(m=2l-n)}}^{2l-3}\frac{1}{n}\frac{1}{m}(n-m)\Big(K^{0}_{n}(s)-K^{0}_{m}(s) \Big)+\frac{2l-2}{2l-1}K^{0}_{2l-1}(s)+\frac{l-2}{2l-2}\Big(K^{0}_{2l-2}(s)-K^{0}_{2}(s)\Big).
\end{split}
\end{equation*}
Any term $(K^{0}_{n}(s)-K^{0}_{m}(s))$ inside the series is such that $m\geq 3$ and $n\geq m+2$. By adding and subtracting equal terms, we obtain
$$
K^{0}_{n}(s)-K^{0}_{m}(s)=\Big(K^{0}_{n}(s)-K^{0}_{n-1}(s)\Big)+\Big(K^{0}_{n-1}(s)-K^{0}_{n-2}(s)\Big)\dots +\Big(K^{0}_{m+1}(s)-K^{0}_{m}(s)\Big), \quad {\rm if}\ m>3
$$
or
$$
K^{0}_{n}(s)-K^{0}_{m}(s)=\Big(K^{0}_{n}(s)-K^{0}_{n-1}(s)\Big)+\Big(K^{0}_{n-1}(s)-K^{0}_{n-2}(s)\Big)\dots +\Big(K^{0}_{5}(s)-K^{0}_{3}(s)\Big), \quad {\rm if}\ m=3.
$$
In both the cases, for point $iii)$  of proposition \ref{lem:propKn} and for remark \ref{rmk:K3}, we infer
$$
K^{0}_{n}(s)-K^{0}_{m}(s)\geq K^{0}_{n}(1)-K^{0}_{m}(1).
$$
Again, for point $i)$ of proposition \ref{lem:propKn}, 
$$
K^{0}_{2l-2}(s)\geq K^{0}_{2l-2}(1)
$$
 and, for lemma \ref{lem:gn}, 
 $$
 \frac{2l-2}{2l-1}K^{0}_{2l-1}(s)-\frac{l-2}{2l-2}K^{0}_{2}(s)
\geq \frac{2l-2}{2l-1}K^{0}_{2l-1}(1)-\frac{l-2}{2l-2}K^{0}_{2}(1)
$$ for any $l\geq 1$.

Therefore, collecting all these bounds, we conclude 
$$
g(s)\geq g(1)=\sum_{\substack{n\geq 2, m\geq 1\\n+m=p-1}}\frac{1}{n}\frac{1}{m}(n-m)K^{0}_{n}(1)=\sum_{\substack{n\geq 2, m\geq 1\\n+m=p-1}}\frac{1}{n}\frac{1}{m}(n-m)\left[\frac{2n}{n+1}-1\right].
$$
The case $p$ even is completely equivalent. Indeed setting $p=2l+2$ we obtain
\begin{equation*}
\begin{split}
g(s)=&\sum_{\substack{n=l+1\\(m=2l+1-n)}}^{2l-2}\frac{1}{n}\frac{1}{m}(n-m)\Big(K^{0}_{n}(s)-K^{0}_{m}s) \Big)\\
&+\frac{2l-3}{2(2l-1)}K^{0}_{2l-1}(s)+\frac{2l-1}{2l}K^{0}_{2l}(s)-\frac{2l-3}{2(2l-1)}K^{0}_{2}(s).
\end{split}
\end{equation*}
Note that when $m=3$ then $n=2l-2=p-4\geq 6$, therefore the term $K_{4}(s)-K_{3}(s)$ is not present.
Arguing as before, we conclude that $g(s)\geq g(1)$.

\qed

Define
$$
S(p)\bydef\sum_{\substack{n\geq 2, m\geq 1\\n+m=p-1}}\frac{1}{n}\frac{1}{m}(n-m)\frac{1}{n+1}.
$$
\begin{lemma}\label{lem:Spneg}
$S(p)$ is decreasing for any $p>4$ and $S(p)<0$ for any $p\geq 11$.
\end{lemma}
\proof
$$
S(p)=\sum_{n=2}^{p-2}\left(\frac{1}{p-1-n} -\frac{1}{n}\right)\frac{1}{n+1}.
$$
Then
 \begin{equation*}
\begin{split}
S(p)-S(p+1)&=\sum_{n=2}^{p-2}\left(\frac{1}{p-1-n} -\frac{1}{n}\right)\frac{1}{n+1}-\sum_{n=2}^{p-1}\left(\frac{1}{p-n} -\frac{1}{n}\right)\frac{1}{n+1}\\
&=\sum_{n=2}^{p-2}\left(\frac{1}{p-1-n} -\frac{1}{p-n}\right)\frac{1}{n+1}-\left(1-\frac{1}{p-1}\right)\frac{1}{p}.
\end{split}
\end{equation*}
Since $\left(\frac{1}{p-1-n} -\frac{1}{p-n}\right)>0$ for any $n=2\dots, p-2$, it holds
\begin{equation*}
\begin{split}
S(p)-S(p+1)&\geq \frac{1}{p-1}\sum_{n=2}^{p-2}\left(\frac{1}{p-1-n} -\frac{1}{p-n}\right)-\frac{p-2}{p(p-1)}\\
&=\frac{1}{p-1}\left[\sum_{n=3}^{p-1}\frac{1}{p-n} -\sum_{n=2}^{p-2}\frac{1}{p-n} \right]-\frac{p-2}{p(p-1)}\\
&=\frac{1}{p-1}\left(1-\frac{1}{p-2} \right)-\frac{p-2}{p(p-1)}=\frac{p-4}{(p-2)(p-1)p}>0\quad \forall p>4.
\end{split}
\end{equation*}
Then $S(p)$ is decreasing for any $p>4$. In particular $S(11)=-\frac{29}{1260}$.

\qed

We are now in the position to prove 
\begin{theorem}\label{thm:bound}
For any $s\in(0,1)$ and any $p\geq 11$
$$
\sum_{\substack{n\geq 2, m\geq 1\\n+m=p-1}}\frac{1}{n}\frac{1}{m}(n-m)K^{0}_{n}(s)>\frac{p-3}{p-2}.
$$
\end{theorem}
\proof
For lemma \ref{lem:bound1}, for any $p\geq 10$
\begin{equation*}
\begin{split}
\sum_{\substack{n\geq 2, m\geq 1\\n+m=p-1}}\frac{1}{n}\frac{1}{m}(n-m)K^{0}_{n}(s)&\geq\sum_{\substack{n\geq 2, m\geq 1\\n+m=p-1}}\frac{1}{n}\frac{1}{m}(n-m)\frac{n-1}{n+1}=\sum_{\substack{n\geq 2, m\geq 1\\n+m=p-1}}\frac{1}{n}\frac{1}{m}(n-m)\left(1-\frac{2}{n+1}\right)\\
&=\sum_{\substack{n\geq 2, m\geq 1\\n+m=p-1}}\left( \frac{1}{m}-\frac{1}{n}\right)-2\sum_{\substack{n\geq 2, m\geq 1\\n+m=p-1}}\frac{1}{n}\frac{1}{m}(n-m)\frac{1}{n+1}.
\end{split}
\end{equation*}
The first sum gives
$$
\sum_{\substack{n\geq 2, m\geq 1\\n+m=p-1}}\left( \frac{1}{m}-\frac{1}{n}\right)=\sum_{n=2}^{p-2}\left(\frac{1}{p-1-n}-\frac{1}{n}\right)=1-\frac{1}{p-2}+\sum_{n=2}^{p-3}\frac{1}{p-1-n}-\sum_{n=2}^{p-3}\frac{1}{n}=\frac{p-3}{p-2}
$$
then, for lemma \ref{lem:Spneg},  the thesis follows.
\qed

\begin{corollary}\label{cor:infinite}
For any $s\in(0,\frac{1}{2}]$ and $p\geq 11$
$$
Q_{0}(p,s)\geq 0.
$$
\end{corollary}
\proof
From \eqref{eq:Qps}

\begin{equation}
\begin{split}
Q_{0}(p,s)=&\left[\sum_{\substack{n\geq 2, m\geq 1\\n+m=p-1}}\frac{1}{n}\frac{1}{m}(n-m)K^{0}_{n}(s)-\frac{p-3}{p-2}\right]\Big(1+\E_0(1-s,q)\Big)^{2}\\
+&\left[\sum_{\substack{n\geq 2, m\geq 1\\n+m=p-1}}\frac{1}{n}\frac{1}{m}(n-m)K^{0}_{n}(1-s)-\frac{p-3}{p-2}\right]\Big(1+\E_0(s,q)\Big)^{2}.\\
\end{split}
\end{equation}
For theorem \ref{thm:bound} both the terms are non negative whenever $p\geq 11$.
\qed

\subsection{Rigorous bound for the finite sum}\label{sec:finite}
 We are now concerning the finite part in the splitting \eqref{eq:separation}.
We wish to prove that
$$
\mathcal I_{0}^{\boldsymbol f}(s,q)=\sum_{p= 4}^{10}Q_{0}(p,s)q^{p-4}\geq 0, \quad \forall s\in(0,1/2), q\in(0,0.9].
$$
By inserting $Q_{0}(p,s)$ in the form \eqref{eq:Qps0}, it results
\begin{equation}
\begin{split}
\mathcal I_{0}^{\boldsymbol f}(s,q)&=\sum_{p=4}^{10}\left[\sum_{\substack{n\geq 2, m\geq 1\\n+m=p}}\frac{1}{n}\frac{1}{m}2(n-m-1)A_{n}(s)-\sum_{\substack{n\geq 2, m\geq 1\\n+m=p-1}}\frac{1}{n}\frac{1}{m}(n-m)A_{n}(s)\right]q^{p-4}\Big(1+\E_0(1-s,q)\Big)^{2}\\
&+\sum_{p=4}^{10}\left[\sum_{\substack{n\geq 2, m\geq 1\\n+m=p}}\frac{1}{n}\frac{1}{m}2(n-m-1)A_{n}(1-s)-\sum_{\substack{n\geq 2, m\geq 1\\n+m=p-1}}\frac{1}{n}\frac{1}{m}(n-m)A_{n}(1-s)\right]q^{p-4}\Big(1+\E_0(s,q)\Big)^{2}\\
&=\mathcal R(s,q)\Big(1+\E_0(1-s,q)\Big)^{2}+\mathcal R(1-s,q)\Big(1+\E_0(s,q)\Big)^{2}
\end{split}
\end{equation}
where  
$$
\mathcal R(s,q)\bydef\sum_{p=4}^{10}\left[\sum_{\substack{n\geq 2, m\geq 1\\n+m=p}}\frac{1}{n}\frac{1}{m}2(n-m-1)A_{n}(s)-\sum_{\substack{n\geq 2, m\geq 1\\n+m=p-1}}\frac{1}{n}\frac{1}{m}(n-m)A_{n}(s)\right]q^{p-4}.
$$
If we perform all the computations, we obtain
\begin{equation*}
\begin{split}
\mathcal R(s,q)=&
-\frac{2}{3}\,s+{\frac {1}{3}}+ \left( 1-\frac{7}{3}\,s+{s}^{2} \right) q+ \left( {
\frac {349}{180}}-{\frac {239}{45}}\,s-\frac{6}{5}\,{s}^{3}+{\frac {43}{10}}\,
{s}^{2} \right) {q}^{2}\\
&+ \left( -{\frac {94}{15}}\,{s}^{3}+{\frac {31}
{10}}+\frac{4}{3}\,{s}^{4}+{\frac {689}{60}}\,{s}^{2}-{\frac {197}{20}}\,s
 \right) {q}^{3}\\
 &+ \left( {\frac {249}{56}}-{\frac {4523}{280}}\,s-{
\frac {10}{7}}\,{s}^{5}+{\frac {4093}{168}}\,{s}^{2}-{\frac {823}{42}}
\,{s}^{3}+{\frac {173}{21}}\,{s}^{4} \right) {q}^{4}\\
&+ \left( \frac{3}{2}\,{s}^
{6}+{\frac {3739}{126}}\,{s}^{4}-{\frac {15367}{630}}\,s+{\frac {3749}
{630}}-{\frac {29941}{630}}\,{s}^{3}-{\frac {143}{14}}\,{s}^{5}+{
\frac {28313}{630}}\,{s}^{2} \right) {q}^{5}\\
&+ \left({\frac {95663}{12600}} -{\frac {218681}{
6300}}\,s+{\frac {474889}{6300}}\,{s}^{2}-{
\frac {620681}{6300}}\,{s}^{3}+{\frac {5125}{63}}\,{s}^{4}-{\frac {
10517}{252}}\,{s}^{5}+{\frac {439}{36}}\,{s}^{6}-{\frac {14}{9}}\,{s}^
{7} \right) {q}^{6}.
\end{split}
\end{equation*}
Note that the coefficient of $q^{0}$ is odd in $s=1/2$.

\begin{theorem}\label{thm:009}
For any $s\in(0,1/2)$ and $q\in(0,0.9]$
$$
\mathcal I_{0}^{\boldsymbol f}(s,q)>0.
$$
\end{theorem}
\proof

For any $q\in(0,1)$ the function $\E(s,q)$ is decreasing in $s$, that is $\E(1,q)\leq \E(s,q)\leq \E(0,q)$. We know that $\mathcal R(s,q)\geq 0 $ for any $q\in[0,1]$ and $s\in[0,\frac{1}{2}]$ while it could be negative for larger values of $s$. Therefore, $\mathcal I_{0}^{\boldsymbol f}(s,q)>0$ for those $s$ where $\mathcal R(1-s)$ is positive, otherwise 
\begin{equation}\label{eq:B}
\mathcal I_{0}^{\boldsymbol f}(s,q)\geq \mathcal B(s,q)\bydef \mathcal R(s,q)\Big(1+\E(1,q)\Big)^{2}+\mathcal R(1-s,q)\Big(1+\E(0,q)\Big)^{2}.
\end{equation}
Hence, it is enough to prove that $ \mathcal B(s,q)>0$. Since $\mathcal B(s,q)$ is made up a finite number of terms and $s,q$ are within  bounded intervals, the problem  $ \mathcal B(s,q)>0$ is, at least theoretically,  well suited to be solved by rigorous numerics. 
 Let us first show that  $ \mathcal B(s,q)>0$ for any $s\in(0,1/2)$ and $q\in(0,0.1]$ then we address the case $q\in[0.1,0.9]$. The reason for this further splitting is that   $\mR(s,q)+\mR(1-s,q)=\mathcal O(q)$, meaning  $\mathcal B(s,q)\to 0$ for any $s$ when $q\to 0$.   This behaviour underlie an obstacle for a rigorous computational scheme to be completely successful: indeed any numerical  tool  based on interval arithmetics will definitely result a non-positive output when computing $\mathcal B(s,q)$ at small values of $q$. Therefore, the case when  $q$ is small deserves a particular treatment.

\vspace{10pt}

\noindent {\bf First computation: $q\in[0,1/10]$}

Let us prove that $\mathcal B(s,q)\geq 0$ for any $s\in[0,1]$, $q\in[0,\frac{1}{10}]$. We compute

\begin{equation*}
\begin{split}
\E(1,q)&\bydef\lim_{s\to1}\E(s,q)=
\frac{2}{q}-1+\frac{2-q}{\ln(1-q)}\\
\E(0,q)&\bydef\lim_{s\to 0}\E(s,q)=
1-\frac{2}{q}-\frac{2-q}{(1-q)\ln(1-q)}.
\end{split}
\end{equation*}
Using the relation $\ln(1-x)^{-1}>-\left(q+\frac{q^{2}}{2}+\frac{q^{3}}{3}+\frac{q^{4}}{4}+\frac{q^{5}}{5} \right)^{-1}$
it follows
$$
\E(1,q)>1-\frac{1}{3}q+\frac{1}{90}q^{3}-\frac{29}{90}q^{4},\quad \forall q\in(0,1).
$$
Similarly,
\begin{equation*}
\begin{split}
\E(0,q)&<1+\frac{1}{3}q+\frac{1}{3}q^{2}+\frac{1}{3}\,{\frac { \left( 58+87\,q+8\,{q}^{2}+15\,{q}^{3}+12\,{q}^{
4} \right) }{ \left( 1-q \right)  \left( 60+30\,q+20\,{q}^{2}+15\,{q}^
{3}+12\,{q}^{4} \right) }}q^{3}\\
&<1+\frac{1}{3}q+\frac{1}{3}q^{2}+0.4q^{3},\quad \forall q\in[0,0.1].
\end{split}
\end{equation*}
By inserting the above estimates into \eqref{eq:B} 
$$
\mathcal B(s,q)\geq \mathcal R(s,q)\Big( 2-\frac{1}{3}q+\frac{1}{90}q^{3}-\frac{29}{90}q^{4}\Big)^{2}+\mathcal R(1-s,q)\Big(2+\frac{1}{3}q+\frac{1}{3}q^{2}+0.4q^{3}\Big)^{2},\quad \forall q\in[0,0.1].
$$
As expected,  the $0^{th}$ order term in $q$ vanishes. Let  be introduced  $\tilde{\mathcal R}(s,q)$, $\tilde {\E}(1,q)$, $\tilde {\E}(0,q)$ so that
$$
\mR(s,q)=\left(\frac{1}{3}-\frac{2}{3}s\right)+q\tilde{\mR}(s,q)
$$
$$
\Big( 2-\frac{1}{3}q+\frac{1}{90}q^{3}-\frac{29}{90}q^{4}\Big)^{2}=4+q\tilde{\E}(1,q),\quad \Big(2+\frac{1}{3}q+\frac{1}{3}q^{2}+0.4q^{3}\Big)^{2}=4+q\tilde{\E}(0,q).
$$
It results
$$
\mathcal B(s,q)\geq q\left[\begin{array}{l}
4\Big(\tilde {\mR}(s,q)+\tilde{\mR}(1-s,q)\Big)+\left(\frac{1}{3}-\frac{2}{3}s \right)\Big(\tilde {\E}(1,q)-\tilde {\E}(0,q) \Big)\\
+q\Big(\tilde {\mR}(s,q)\tilde {\E}(1,q)+\tilde{\mR}(1-s,q)\tilde {\E}(0,q) \Big)
\end{array} \right].
$$
For the choice $\Delta(s)=0.02$, $\Delta(q)=0.01$, define the intervals
$$
s_{j}\bydef[(j-1)\Delta(s),j\Delta_{s}] ,\quad  j=1,\dots,25,\qquad q_{k}\bydef[(k-1)\Delta_{q},k\Delta_{q}] ,\quad k=1,\dots,10
$$
so that $\bigcup_{j} s_{j}\supset [0,1/2]$ and $\bigcup_{k} q_{k}\supset [0,1/10]$, Then, for any $j,k$ in their range, we compute using  interval arithmetics the quantities
\begin{equation}
\begin{split}
\mathcal M(j,k)=&4\Big(\tilde {\mR}(s_{j},q_{k})+\tilde{\mR}(1-s_{j},q_{k})\Big)+\left(\frac{1}{3}-\frac{2}{3}s_{j} \right)\Big(\tilde {\E}(1,q_{k})-\tilde {\E}(0,q_{k}) \Big)\\
&+q_{k}\Big(\tilde {\mR}(s_{j},q_{k})\tilde {\E}(1,q_{k})+\tilde{\mR}(1-s_{j},q_{k})\tilde {\E}(0,q_{k}) \Big).
\end{split}
\end{equation}
It results $\min_{j,k}(\min(\mathcal M(j,k)))>0.2744$, proving that $\mathcal B(s,q)> 0$ for any $s\in\Big(0,\frac{1}{2}\Big)$, $q\in\Big(0,\frac{1}{10}\Big]$.

\vspace{10pt}
\noindent {\bf Second computation, $q\in[0.1,0.9]$}

We are now concern the remaining part, that is when $q\in[\frac{1}{10},\frac{9}{10}]$. We rigorously compute the lower bound for $\mathcal B(s,q)$, as given in \eqref{eq:B},  for $s\in(0,1)$, $q\in[0.1,0.9]$. With the choice of $\Delta_{s}=2\cdot 10^{-3}$ and $\Delta_{q}=2\cdot 10^{-4}$ it results that $\mathcal B(s,q)\geq 0.0013$.

\qed

We can now state the theorem:
\begin{theorem}
For any value of the energy $E$, the apsidal angle $\Delta_{0}\theta(u)$ is monotonically increasing as function of the angular momentum.
\end{theorem}

\proof 
From proposition \ref{prop:0.9}, corollary \ref{cor:infinite}  and theorem \ref{thm:009} it follows that $\mathcal I_{0}(s,q)>0$ for any $s\in(0,1)$ and $q\in(0,1)$. Then, for lemma \ref{lem:3/2}, the derivative in \eqref{eq:der} is negative for any $q\in(0,1)$. Since $q=q(\ell)$ is decreasing, it follows that 
$$
\frac{d}{d\ell}\Delta_{0}\theta(u)>0, \quad \forall \ell\in(0,\ell_{max}).
$$
\qed

\section{Appendix}

This Appendix is intended to collect all the technical results concerning the function $A_{n}(s)$ and $K^{\al}_{n}(s,q)$.

\subsection{Properties of $A_{n}(s)$}
$$
A_{n}(s)=\frac{1-(1-s)^{n-1}}{s},\quad n\geq 2, s\in(0,1)
$$

\begin{lemma}\label{lem:propA}
\begin{itemize}
\item[]
\item[i)] $A_{2}(s)=1$, $A_{n}(1)=1$ for any $n\geq 2$.
\item[ii)] $A_{n+1}(s)-A_{n}(s)=(1-s)^{n-1},\quad A_{n+1}(1-s)-A_{n}(1-s)=s^{(n-1)}$.
\item[iii)] for any $n\geq 3$ and $s\in(0,1)$, $A_{n}(s)$ is decreasing and convex, i.e. $A_{n}'(s)<0$ and $A_{n}''(s)>0$. 
\item[iv)] for any $n>m$ $A'_{n}(s)<A'_{m}(s)$  for any $s\in(0,1)$.
\end{itemize}
\end{lemma}

\proof
 
$i) , ii)$ obvious.

$iii)$
$$
\frac{d}{ds}A_{n}(s)=\frac{-1+(1-s)^{n-1}+(n-1)(1-s)^{(n-2)}s}{s^{2}}.
$$
Looking at the numerator, say $N(s)$, we see that $N(0)=0$, $N(1)=-1$ and
$$
N'(s)=-(n-1)(1-s)^{n-2}-(n-1)(n-2)(1-s)^{n-3}s+(n-1)(1-s)^{(n-2)}<0
$$
hence $A_{n}(s)$ is decreasing.

$$
\frac{d^{2}}{ds^{2}}A_{n}(s)=\frac{2s-2s(1-s)^{(n-1)}-2(n-1)(1-s)^{n-2}s^{2}-(n-2)(n-1)(1-s)^{n-3}s^{3}}{s^{4}}
$$
$$
=\frac{2-2(1-s)^{(n-1)}-2(n-1)(1-s)^{n-2}s-(n-2)(n-1)(1-s)^{n-3}s^{2}}{s^{3}}.
$$
Looking at the numerator, $N(s)$, we note that $N(0)=0, N(1)=2$ and
$$
N'=(n-3)(n-2)(n-1)(1-s)^{n-3}>0, s\in(0,1)
$$
hence $A_{n}(s)$ is convex.

$iv)$
\begin{equation*}
\begin{split}
\frac{d}{dn}A'_{n}(s)&=\frac{(1-s)^{n-1}\log(1-s)+(1-s)^{n-2}s+(n-1)(1-s)^{n-2}s\log(1-s)}{s^{2}}\\
&=\frac{(1-s)^{(n-2)}}{s^{2}}(\log(1-s)(1-s)+s+(n-1)s\log(1-s))\\
&=\frac{(1-s)^{(n-2)}}{s^{2}}(\log(1-s)+s+(n-2)s\log(1-s))<0.
\end{split}
\end{equation*}

\qed

\subsection{Properties of $K^{\al}_{n}(s)$}

$$
K^{\al}_{n}(s)=2\frac{n-\al}{n+1}A_{n+1}(s)-A_{n}(s), \quad s\in(0,1), \al\in[0,1), n\geq 2.
$$

\begin{lemma}\label{lem:propKn}
For any $\al\in[0,1)$ and $n\geq 2$:
\begin{itemize}
\item[i)] $K^{\al}_{n}(s)$ is decreasing and convex  for any $s\in(0,1)$.
\item[ii)] $K^{\al}_{n+1}(s)-K^{\al}_{n}(s)>0$ for any $s\in(0,1)$.
\end{itemize}
For any $\al\in[0,1)$ and $n\geq 4$
\begin{itemize}
\item[iii)] $K^{\al}_{n+1}(s)-K^{\al}_{n}(s)$ is decreasing for any $s\in(0,1)$.

\end{itemize}
\end{lemma}

\proof

$i)$
By direct computation the property is proved for $n=2,3$.
Then, 
\begin{equation}
\begin{split}
K^{\al}_{n}(s)=&\frac{2(n-\al)}{n+1}A_{n+1}(s)-\frac{1-(1-s)^{n-1}+(1-s)^{n}-(1-s)^{n}}{s}\\
&\left(\frac{2(n-\al)}{n+1}-1 \right)A_{n+1}(s)+(1-s)^{n-1},
\end{split}
\end{equation}
hence
$$
\frac{d}{ds}K^{\al}_{n}(s)=\left(\frac{n-2\al-1}{n+1} \right)A'_{n+1}(s)-(n-1)(1-s)^{n-2}<0,\quad  \forall s\in(0,1), n\geq 4.
$$ 
Also,
$$
\frac{d^{2}}{ds^2}K^{\al}_{n}(s)=\left(\frac{n-2\al-1}{n+1} \right)A''_{n+1}(s)+(n-2)(n-1)(1-s)^{n-3}>0,\quad  \forall s\in(0,1), n\geq 4.
$$

$ii)$
The case $n=2$ gives  $K^{\al}_{3}(s)-K^{\al}_{2}(s)=\frac{5}{6}-{\frac {13}{6}}\,s+\frac{3}{2}\,{s}^{2}-\frac{1}{6}\,a+\frac{5}{6}\,\al s-\frac{1}{2}\al {s}^{2}$ that is positive for any $s\in(0,1)$ and $\al\in[0,1)$.

In the general case 
$$
 K^{\al}_{n+1}(s)-K^{\al}_{n}(s)=
 \frac{2(n+1-\al)}{n+2}A_{n+2}(s)-\frac{3n+1-2\al}{n+1}A_{n+1}(s)+A_{n}(s)
 $$
 $$
 =\frac{1}{(n+1)(n+2)}\left[2(n+1)(n+1-\al)A_{n+2}(s)-(3n+1-2\al)(n+2)A_{n+1}(s)+(n+1)(n+2)A_{n}(s) \right].
 $$
 Using twice property $ii)$ of lemma \ref{lem:propA} we obtain
 $$
 K^{\al}_{n+1}(s)-K^{\al}_{n}(s)=\frac{1}{(n+1)(n+2)}\cdot
 $$
 $$
 \left[2(\al+1)A_{n}(s)+2(n+1)(n+1-\al)(1-s)^{n}-(n^{2}+3n-2\al)(1-s)^{n-1} \right].
 $$
 Denote by $\eta^{\al}_{n}(s)\bydef2(n+1)(n+1-\al)(1-s)^{n}-(n^{2}+3n-2\al)(1-s)^{n-1}$. Its derivative  $\dot \eta^{\al}_{n}(s)=(1-s)^{n-2}((n^{2}+3n-2\al)(n-1)-2n(n+1)(n+1-\al)(1-s))$ is zero for $(1-s)=\frac{(n^{2}+3n-2\al)(n-1)}{2n(n+1)(n+1-\al)}$ where it corresponds the minimum for $\eta^{\al}_{n}(s)$. We infer
 \begin{equation}
\begin{split}
 \eta^{\al}_{n}(s)&\geq 2(n+1)(n+1-\al)\frac{(n^{2}+3n-2\al)^{n}(n-1)^{n}}{(2n)^{n}(n+1-\al)^{n}(n+1)^{n}}-\frac{(n^{2+3n-2\al})^{n}(n-1)^{n-1}}{(2n)^{n-1}(n+1-\al)^{n-1}(n+1)^{n-1}}\\&=-\frac{(n^{2}+3n-2\al)^{n}(n-1)^{n-1}}{2^{n-1}n^{n}(n+1-\al)^{n-1}(n+1)^{n-1}}.
 \end{split}
 \end{equation}
 For any $\al\in[0,1)$ we can uniformly  bound 
 $$
 \eta_{n}^{\al}(s)\geq \eta(n)\bydef-\frac{(n^{2}+3n)^{n}(n-1)^{n-1}}{2^{n-1}n^{n}(n)^{n-1}(n+1)^{n-1}}=-\frac{(n+3)^{n}(n-1)^{n-1}}{(2n)^{n-1}(n+1)^{n-1}}.
 $$
 The last ratio is monotonically decreasing for $n\geq 1$, and it equals $\frac{3}{2}$ at $n=3$. Therefore $\eta^{\al}_{n}(s)\geq -\frac{3}{2}$ for any $s\in(0,1)$ and $n\geq 3$, $\al\in[0,1)$. Since $A_{n}(s)\geq 1$ for any $n$ and $s\in(0,1)$, it follows 
 $$
 K^{\al}_{n+1}(s)-K^{\al}_{n}(s)\geq\frac{1}{(n+1)(n+2)}\left(2-\frac{3}{2}\right)>0.
 $$

$iii)$
The case  $n=4,5$ is proved by direct computation. 
%
 For $n\geq 6$ 
 $$
 \frac{d}{ds}\Big(K^{\al}_{n+1}(s)-K^{\al}_{n}(s)\Big)=\frac{1}{(n+1)(n+2)}\cdot
 $$
 $$
 \left[(2+\al)A'_{n}(s)-2n(n+1)(n+1-\al)(1-s)^{n-1}+(n-1)(n^{2+}3n-2\al)(1-s)^{n-2} \right].
 $$
 Denote by $\nu^{\al}_{n}(s)\bydef-2n(n+1)(n+1-\al)(1-s)^{n-1}+(n-1)(n^{2+}3n-2\al)(1-s)^{n-2} $. Then, arguing as before, we prove that 
 $$
 \max_{s\in(0,1)}(\nu^{\al}_{n}(s))=\frac{(n-2)^{n-2}(n^{2}+3n-2\al)^{n-1}}{(2n)^{n-2}(n+1)^{n-2}(n+1-\al)^{n-2}}.
 $$
 Therefore, for any $\al\in[0,1)$,
 $$
 \nu^{\al}_{n}(s)\leq \nu(n)\bydef \frac{(n-2)^{n-2}(n^{2}+3n)^{n-1}}{(2n)^{n-2}(n+1)^{n-2}(n)^{n-2}}= \frac{(n-2)^{n-2}n(n+3)^{n-1}}{(2n)^{n-2}(n+1)^{n-2}}.
 $$ 
 The function $\nu({n})$ is decreasing in $n$ and $\nu(6)<2$.
 Since $A'_{n}(s)\leq -1$ for any $s$ and $n$, it follows that $K^{\al}_{n+1}(s)-K^{\al}_{n}(s)$ is decreasing in $s$ for any $n\geq 6$.
 
\qed

\begin{remark}\label{rmk:K3}
$K^{\al}_{4}(s)-K^{\al}_{3}(s)$ is not decreasing in $s\in(0,1)$ and the minimum is not achieved at $s=1$. However the function $K^{0}_{5}(s)-K^{0}_{3}(s)$, for $s\in [0,1]$ attains the minimum at $s=1$. Indeed
$$
K^{0}_{5}(s)-K^{0}_{3}(s)={\frac {11}{6}}-{\frac {43}{6}}\,s+{\frac {67}{6}}\,{s}^{2}-{\frac {22}{3}}\,{s}^{3}+\frac{5}{3}\,{s}^{4}
$$  
and the function has a local minimum and a local maximum at $s_{1}$, $s_{2}$ respectively, where
$s_{1}\in I_{1}=[0.7230,0.7233]$ and $s_{2}\in I_{2}=[  0.871,0.872]$. For any $s\in I_{1}$, $K^{0}_{5}(s)-K^{0}_{3}(s)\in [  0.1672,    0.1785] >\frac{1}{6}=K^{0}_{5}(1)-K^{0}_{3}(1)$.
\end{remark}

\begin{lemma}\label{lem:gn}
For any $n\geq 2$, the function
$$
g_{n}(s)\bydef\frac{n-1}{n}K^{0}_{n}(s)-\frac{n-3}{2(n-1)}K^{0}_{2}(s)
$$
is monotonically decreasing for any $s\in[0,1]$.
\end{lemma}
\proof

$$
g_{n}(s)=\frac{n-1}{n}K^{0}_{n}(s)-\frac{n-3}{2(n-1)}\left(\frac{4}{3}A_{3}(s)-1 \right)
$$ 
then, since $A'_{3}(s)=-1$,
$$
\dot g_{n}(s)=\frac{n-1}{n}\frac{d}{ds}K^{0}_{n}(s)+\frac{4}{3}\frac{n-3}{2(n-1)}.
$$
For point $i)$ of lemma \ref{lem:propKn}, for any $s\in(0,1)$,
$$
\dot g_{n}(s)\leq -\frac{n-1}{n}\left(\frac{n-1}{n+1}\right)+\frac{2(n-3)}{3(n-1)}=\frac{-3(n-1)^{3}+2n(n-3)(n+1)}{3n(n+1)(n-1)}<0, \quad \forall n\geq 2.
$$
\qed

\bibliographystyle{unsrt}

\end{document}